\documentclass{article}
\usepackage[margin=1.25in]{geometry}   

\author{John K. Sikora \\ sikorajk@corning.com}
\title{On Calculating the Coefficients of a Polynomial Generated Sequence Using the \wnts}

\usepackage{hyperref}   
\usepackage{xspace}     
\usepackage{amsthm}     
\usepackage{amsmath}    
\usepackage{amsfonts}    
\usepackage{twoopt}     
\usepackage{verbatim}   
\usepackage{mathtools}  
\usepackage{diagbox}    
\usepackage{cleveref}   

\newcommand{\wns}{Worpitzky Numbers\xspace}
\newcommand{\snsk}{Stirling Number of the Second Kind\xspace}
\newcommand{\snskm}{\ifmmode S(n,k) \else $S(n,k)$\xspace \fi} 
\newcommand{\mwnt}{Mirrored Worpitzky Number Triangle\xspace}
\newcommand{\wnt}{Worpitzky Number Triangle\xspace}
\newcommandtwoopt{\mwntm}[2][n][k]{\ifmmode MWNT({#1},{#2}) \else $MWNT({#1},{#2})$\xspace \fi} 
\newcommandtwoopt{\awntm}[2][n][k]{\ifmmode AWNT({#1},{#2}) \else $AWNT({#1},{#2})$\xspace \fi} 
\newcommand{\wnts}{Worpitzky Number Triangles\xspace}
\newcommand{\awnt}{Alternate Worpitzky Number Triangle\xspace}
\newcommand{\dt}{difference triangle\xspace}
\newcommand{\dtc}{Difference Triangle\xspace}
\newcommand{\ps}{polynomial sequence\xspace}
\newcommand{\psc}{Polynomial Sequence\xspace}
\newcommand{\pspl}{polynomial sequences\xspace}
\newcommand{\mdc}{Main Diagonal\xspace}
\newcommand{\md}{main diagonal\xspace}
\newcommand{\efdt}{Euler's Finite Difference Theorem\xspace}
\newcommand{\gaq}{Gould and Quaintance, (2016)\xspace}
\newcommand{\oeis}{OEIS\xspace}
\newcommand{\mwnta}{A028246\xspace}
\newcommand{\awnta}{A019538\xspace}
\newcommand{\wnta}{A130850\xspace}
\newcommand{\snska}{A008277\xspace}
\newcommand{\moem}{\hspace{-1em}}  
\newcommand{\mxem}{\hspace{-2.95em}}

\let\oldFootnote\footnote
\newcommand\nextToken\relax

\renewcommand\footnote[1]{%
    \oldFootnote{#1}\futurelet\nextToken\isFootnote}

\newcommand\isFootnote{%
    \ifx\footnote\nextToken\textsuperscript{,}\fi}

\newtheoremstyle{jsdefn}
{10pt}
{10pt}
{}
{}
{\bfseries}
{:}
{1em}
{}

\newtheoremstyle{jsexample}
{14pt}
{14pt}
{}
{}
{\bfseries \itshape}
{  --}
{1em}
{}

\theoremstyle{jsdefn}
\newtheorem{defn}{Definition} 
\theoremstyle{jsexample}
\newtheorem{examp}{Example}
\theoremstyle{plain}

\begin{document}
\numberwithin{equation}{section}

\maketitle 

\begin{abstract}
In this paper we show that two of the \wnts, \oeis \mwnta and \awnta, may each be used in look-up table fashion, along with specific diagonals of a \ps's \dt to easily solve for the unknown coefficients of the sequence. This is accomplished by using a method that isolates each of the coefficients as a single unknown in a series of simple linear equations. The method is first applied to a sequence generated using integer indexes with a starting index of 0, using the \awnt, \awnta. Although the numbers in \awnta are less commonly referred to as a \wnt, a justification for such a reference is briefly presented. Next, the method is applied to a sequence generated using integer indexes with a starting index of 1, using the \mwnt, \mwnta.  The method is then extended to solve for the coefficients of a sequence generated with input data consisting of an arbitrary starting number and an arbitrary constant differential.
\end{abstract}

\section{Introduction} \label{sec:intro}
It is well known that the first step in using what is perhaps the most popular method to find the unknown generating polynomial of a sequence is to calculate the \dt of the sequence by subtracting successive terms of the sequence. The next row in the \dt is formed in a like manner, by subtracting successive terms. If the $d^{th}$ row of the \dt, as defined below in \Cref{def:ddt}, stays constant for a sufficient number of terms, then it is known that that the generating sequence is a polynomial of degree $d$. 

Once it is established that the sequence is generated via a polynomial of degree $d$, there are various methods that may be used to obtain a formula for the generating polynomial. Such methods include generating and solving linear equations for the coefficients of each power in the polynomial, or using Newton's Divided Difference Formula.

We present what we believe is an easy and newly described approach of formulating the unknown polynomial using a \wnt and the \md of the \dt of the sequence, given, or assuming that the sequence is generated with input data consisting of integers starting at either 0 or 1. We then show how the method may be generalized to allow it to be used on a sequence with input data consisting of an arbitrary starting number and an arbitrary constant differential.

\section{Outline} \label{sec:outl}
The remainder of this paper consists of a number of sections. Following this outline, in the next section the definitions, notations, and the specific formulas that are used in the paper are given. In addition, references to alternate versions of the \wnts are presented. 

Next, three examples are given, with only the practical calculations shown. We feel that this is desirable in order to show the ease of using the method. The three examples show the calculations for a \ps generated with input data consisting of integers starting at 0, then for a sequence generated with input data consisting of integers starting at 1, and finally for a sequence generated with input data starting at 3.3 with an increment of 0.1. 

After that, in the next section we provide the mathematical basis for the method, and in the following two sections we first present the the full rendition of one of the examples, and then a partial rendition of another example. We feel that this will make the mathematics behind the method more apparent. In the final section, we present our closing remarks. 

\section{Definitions, Notation, and Existing Terminology} \label{sec:term}

\begin{defn} \label{def:dmwnt}    
\textbf{\mwnt} or \boldmath \mwntm \unboldmath -- The triangle formed from the numbers, $n$ and $k$, in \oeis \cite{bboeis} \mwnta \cite{bbmwnta} as shown in \cite[Example]{bbmwnta} and in \Cref{tbl:mwntt}. One formula for the numbers in the \mwnt is $(k-1)! \cdot \snskm$, where $S(n,k)$ is the \snsk \cite{bbqgsn} \cite{bbcmsn}. A triangle of these numbers is given in \oeis \snska \cite{bbsnsk}. 

\vspace{0.1em}
The specific formula for \mwntm used in this paper, equivalent to the one given above, is:

\begin{equation} \label{eq:mwnt}
\mwntm = \frac{1}{k} \sum_{i=0}^{k} (-1)^{k-i} \binom{k}{i} i^n \ ; \ n \geq 1, \ k \geq 1
\end{equation}

\vspace{0.1em} 
The mirror image of \mwnta was recently referred to as the \wnt by Vandervelde \cite{bbvdv}, and we yield to that reference by using the term ``Mirrored'' in our definition. The referenced triangle may be found on the \oeis as \wnta \cite{bbmiwnt}. However, it should be noted that in \oeis \mwnta, \wnta is referred to as ``The mirror image of the Worpitzky triangle'' \cite[Comments]{bbmwnta}. 

\vspace{0.1em} 
In addition, what we refer to as the \mwnt (\mwnta) appears elsewhere in the \oeis, such as in OEIS sequence A005460 \cite[Links]{bbrsa} (\cite{bbrsanps} provides a direct link). A005460 is described \cite[Comments]{bbrsa} as: ``third external diagonal of Worpitzky triangle \mwnta''.\footnote{Although we had figured out that the first diagonal in what we would eventually call the \mwnt, \mwntm, is $(n-1)!$, and that the second diagonal is $n!/2$, we were perplexed about the third diagonal which is 1, 7, 50, 390, 3360, etc. A search on the OEIS turned up A005460, which referenced \mwnta, and all of the succeeding diagonals that we checked matched. Since these numbers were readily available on the \oeis in look up table form, we decided to write this paper.} Obviously, the use of the term \wnt (or similar) varies.
\end{defn}

\begin{table}[!htbp]
\centering
\caption{The \mwnt, \oeis \mwnta, with zeros for $n < k$}
\label{tbl:mwntt}
\begin{tabular}{l|rrrrrrrrr} 
\backslashbox{n}{\vspace{-1.5em} k} & \ 1 & 2 & 3 & 4 & 5 & 6 & 7 & 8 & 9 \\ \hline
1 & 1 & 0   & 0    & 0     & 0      & 0      & 0      & 0      & 0 \\
2 & 1 & 1   & 0    & 0     & 0      & 0      & 0      & 0      & 0 \\
3 & 1 & 3   & 2    & 0     & 0      & 0      & 0      & 0      & 0 \\
4 & 1 & 7   & 12   & 6     & 0      & 0      & 0      & 0      & 0 \\
5 & 1 & 15  & 50   & 60    & 24     & 0      & 0      & 0      & 0 \\
6 & 1 & 31  & 180  & 390   & 360    & 120    & 0      & 0      & 0 \\
7 & 1 & 63  & 602  & 2100  & 3360   & 2520   & 720    & 0      & 0 \\
8 & 1 & 127 & 1932 & 10206 & 25200  & 31920  & 20160  & 5040   & 0 \\
9 & 1 & 255 & 6050 & 46620 & 166824 & 317520 & 332640 & 181440 & 40320 \\
\end{tabular}
\end{table}

\begin{defn} \label{def:dawnt}    
\textbf{\awnt} or \boldmath \awntm \unboldmath -- The triangle formed from the numbers, $n$ and $k$, in \oeis \awnta \cite{bbawnta} as shown in \cite[Example]{bbawnta} and in \Cref{tbl:awntt}. Perhaps providing justification for referring to this triangle as a \wnt comes from \gaq \cite[Equation 11.3]{bbqg}, who provide an equation for \wns in general. A specific case is mentioned \cite{bbqgsc} which results in the numbers in the \awnt, with a formula given as as $k! \cdot \snskm$.

\vspace{0.1em} 
The specific formula for \awntm used in this paper, equivalent to the one given above, is:

\begin{equation} \label{eq:awnt}
\awntm = \sum_{i=0}^{k} (-1)^{k-i} \binom{k}{i} i^n \ ; \ n \geq 1, \ k \geq 1
\end{equation}
\end{defn}

\begin{table}[!htbp]
\centering
\caption{The \awnt, \oeis \awnta, with zeros for $n < k$} 
\label{tbl:awntt}
\begin{tabular}{l|rrrrrrrrr} 
\backslashbox{n}{\vspace{-1.5em} k} & \ 1 & 2 & 3 & 4 & 5 & 6 & 7 & 8 & 9 \\ \hline
1 & 1 & 0   & 0     & 0      & 0      & 0       & 0       & 0       & 0 \\
2 & 1 & 2   & 0     & 0      & 0      & 0       & 0       & 0       & 0 \\
3 & 1 & 6   & 6     & 0      & 0      & 0       & 0       & 0       & 0 \\
4 & 1 & 14  & 36    & 24     & 0      & 0       & 0       & 0       & 0 \\
5 & 1 & 30  & 150   & 240    & 120    & 0       & 0       & 0       & 0 \\
6 & 1 & 62  & 540   & 1560   & 1800   & 720     & 0       & 0       & 0 \\
7 & 1 & 126 & 1806  & 8400   & 16800  & 15120   & 5040    & 0       & 0 \\
8 & 1 & 254 & 5796  & 40824  & 126000 & 191520  & 141120  & 40320   & 0 \\
9 & 1 & 510 & 18150 & 186480 & 834120 & 1905120 & 2328480 & 1451520 & 362880 \\
\end{tabular}
\end{table}

\begin{defn} \label{def:ps}    
\textbf{\psc} -- A sequence, $a_i, a_{i+1},a_{i+2}$, etc., generated by the polynomial of finite degree $d$ and written in long form as $c_d x^d + c_{d-1} x^{d-1} + c_{d-2} x^{d-2} + \cdots + c_1 x + c_0$. In this paper the more compact form, $\sum_{j=0}^d c_j x^j$ will primarily be used. The $x$ values may be integers or real numbers with an arbitrary starting value and an arbitrary constant differential.
\end{defn}

\begin{defn} \label{def:ddt}    
\textbf{The \dtc and the \mdc} -- The difference triangle of a sequence is the triangle formed by subtracting the preceding element of a sequence from the current element, and continuing this process for successive rows. An example is shown in \Cref{tbl:gdt}. Note that the row containing the sequence values,  $a_0,a_1,a_2$, etc. is row number 0, and the succeeding rows are numbered $1,2,3,$ etc. The left-most diagonal is shown in bold and is known as the \md.
\end{defn}

\begin{table}[!htbp]
\centering
\caption{The General Difference Table for a Sequence}
\label{tbl:gdt}
\begin{tabular}{cccccccc} 
$\boldsymbol{a_0}$ & & $a_1$ & & $a_2$ & & \hspace{0.5em} $a_3$ \\
& \hspace{1em} $\boldsymbol{a_1-a_0}$ & & $a_2-a_1$ & & \hspace{-0.5em} $a_3-a_2$ \\
& & $\boldsymbol{a_2-2a_1+a_0}$ & & $a_3-2a_2+a_1$ \\
& & & $\boldsymbol{a_3-3a_2+3a_1-a_0}$ \\
\end{tabular}
\end{table}

\begin{defn} \label{def:dztoz}    
\boldmath $0^0 = 1$ \unboldmath -- More specifically, $x^0 = 1$ for all $x$, per Graham, Knuth, and Patashnik, (1994) \cite{bbcm}. This is common when using binomials, and it allows for the $c_0x^0$ term to be $c_0$ when $x=0$ (as is necessary) in the compact formula for the \ps given in \Cref{def:ps}.
\end{defn}

\section{Examples of Using the Method} \label{sec:ex}

The following examples provide a brief description of how to use the method.

\begin{examp} \label{ex:awnt} 
We first consider the sequence generated by the polynomial: \[4x^6 + 5x^5 + 6x^4 + 7x^3 + 8x^2 + 9x + 10; \ x \in 0..7\] The difference table for this sequence is given in \Cref{tbl:awntex}, with the \md in bold. The sequence values are in row 0 in accordance with \Cref{def:ddt}. In this example it is either given or assumed that the the integers starting with 0 were used to generate the sequence. Therefore, we will use the \awnt and the \md values to easily solve for the ``unknown'' coefficients. 

\begin{table}[!htbp]
\centering
\caption{The Difference Table for $4x^6 + 5x^5 + 6x^4 + 7x^3 + 8x^2 + 9x + 10; \ x \in 0..7$}
\label{tbl:awntex}
\begin{tabular}{cccccccccccccccl} 
0 & & \moem 1 & & \moem 2 & & \moem 3 & & \moem 4 & & \moem 5 & & \moem 6 & & \moem 7 & $\leftarrow x$ \\
\textbf{10} & & \moem 49 & & \moem 628 & & \moem 4915 & & \moem 23662 & & \moem 83005 & & \moem 235144 & & \moem 571903 \\
& \moem \textbf{39} & & \moem 579 & & \moem 4287 & & \moem 18747 & & \moem 59343 & & \moem 152139 & & \moem 336759 \\
& & \moem \textbf{540} & & \moem 3708 & & \moem 14460 & & \moem 40596 & & \moem 92796 & & \moem 184620 \\
& & & \moem \textbf{3168} & & \moem 10752 & & \moem 26136 & & \moem 52200 & & \moem 91824 \\
& & & & \moem \textbf{7584} & & \moem 15384 & & \moem 26064 & & \moem 39624 \\
& & & & & \moem \textbf{7800} & & \moem 10680 & & \moem 13560 \\
& & & & & & \moem \textbf{2880} & & \moem 2880 \\
\end{tabular}
\end{table}

\vspace{0.1em} 
Since the integers that generate the sequence start at 0, we know that $c_0=10$, as 10 is the initial value in row 0. Row 6 of the \dt is constant, which, of course, matches the degree, $d$, of the polynomial. Therefore, we start at $n,k=d=6$ in \Cref{tbl:awntt}, and use column $k$ to calculate $c_{n=k}$ of the polynomial using the entries in the table as multipliers for each $c_n, n \in 1..6$. The $c_n$ values are multiplied by the values in row $n$ of \Cref{tbl:awntt} for that column. As can be seen in the calculations below, this process turns the rows of the table into columns, and the columns of the table into rows, when presented in the manner shown.

We will move backwards along the columns, $k$, in turn calculating the $c_{n=k}$ values as we go. The reasons for these steps and the mathematical relationships between the coefficients and the multiplier values in the table will be shown in \Cref{ss:mdawnt}. The multipliers (elements of the table) appear in parenthesis below, with multipliers of 0 not shown. Thus, we have:
\[ \begin{array}{ccccccccccccccccc}    
2880 & = & c_6(720) & & & & & & & & & & & \Rightarrow & c_6 & = & 4 \\
7800 & = & 4(1800) & + & c_5(120) & & & & & & & & & \Rightarrow & c_5 & = & 5 \\
7584 & = & 4(1560) & + & 5(240) & + & c_4(24) & & & & & & & \Rightarrow & c_4 & = & 6 \\
3168 & = & 4(540) & + & 5(150) & + & 6(36) & + & c_3(6) & & & & & \Rightarrow & c_3 & = & 7 \\
540 & = & 4(62) & + & 5(30) & + & 6(14) & + & 7(6) & + & c_2(2) & & & \Rightarrow & c_2 & = & 8 \\
39 & = & 4(1) & + & 5(1) & + & 6(1) & + & 7(1) & + & 8(1) & + & c_1(1) & \Rightarrow & c_1 & = & 9 \\
\end{array} \]
Since we already know that $c_0=10$, we have the complete solution for \Cref{ex:awnt}.
\end{examp}

\begin{examp} \label{ex:mwnt} 
In this example, we show how to directly calculate the coefficients of a \ps given that the integers starting with 1 (instead of 0) are used to generate the sequence. We will use the values in the \mwnt (instead of the \awnt) and in the \md of the \dt.  

\vspace{0.1em} 
We could have repeated \Cref{ex:awnt} using the next diagonal (adjacent to the \md) of \Cref{tbl:awntex}, but we elect to use a different sequence to add more variety to the examples. We consider the sequence generated by the polynomial: \[2x^6 + 3x^5 + 5x^4 + 7x^3 + 11x^2 + 13x + 17; \ x \in 1..8\] 

The difference table for this sequence is given in \Cref{tbl:wntex}, with the \md in bold. Again, the sequence values are in row 0 in accordance with \Cref{def:ddt}.  

\begin{table}[!htbp]
\centering
\caption{The Difference Table for $2x^6 + 3x^5 + 5x^4 + 7x^3 + 11x^2 + 13x + 17; \ x \in 1..8$}     
\label{tbl:wntex}
\begin{tabular}{cccccccccccccccl} 
\moem 1 & & \moem 2 & & \moem 3 & & \moem 4 & & \moem 5 & & \moem 6 & & \moem 7 & & \moem 8 & $\leftarrow x$ \\
\moem \textbf{58} & & \moem 447 & & \moem 2936 & & \moem 13237 & & \moem 44982 & & \moem 125123 & & \moem 300772 & & \moem 647481 \\
& \moem \textbf{389} & & \moem 2489 & & \moem 10301 & & \moem 31745 & & \moem 80141 & & \moem 175649 & & \moem 346709 \\
& & \moem \textbf{2100} & & \moem 7812 & & \moem 21444 & & \moem 48396 & & \moem 95508 & & \moem 171060 \\
& & & \moem \textbf{5712} & & \moem 13632 & & \moem 26952 & & \moem 47112 & & \moem 75552 \\
& & & & \moem \textbf{7920} & & \moem 13320 & & \moem 20160 & & \moem 28440 \\
& & & & & \moem \textbf{5400} & & \moem 6840 & & \moem 8280 \\
& & & & & & \moem \textbf{1440} & & \moem 1440 \\
\end{tabular}
\end{table}

\vspace{0.1em}
Row 6 of the \dt is constant, again matching the degree, $d$, of the polynomial. We start at $n,k=d+1=7$ in \Cref{tbl:mwntt}, and use column $k$ to calculate $c_{n-1=k-1}$ of the polynomial using the entries in the table as multipliers for each $c_{n-1}, n \in 1..7$. The $c_{n-1}$ values are multiplied by the values in row $n$ of \Cref{tbl:mwntt} for that column. Again, we will move backwards along the columns, calculating the $c_{n-1=k-1}$ values as we go. The reasons for these steps and the mathematical relationships between the coefficients and the multiplier values in the table will be shown in \Cref{ss:mdmwnt}.
The multipliers (elements of the table) appear in parenthesis below, with multipliers of 0 not shown. Thus, we have:
\[ \arraycolsep=3.0pt
\begin{array}{ccccccccccccccccccc} 
1440 & = & c_6(720) & & & & & & & & & & & & & \Rightarrow & c_6 & = & 2 \\
5400 & = & 2(2520) & + & c_5(120) & & & & & & & & & & & \Rightarrow & c_5 & = & 3 \\
7920 & = & 2(3360) & + & 3(360) & + & c_4(24) & & & & & & & & & \Rightarrow & c_4 & = & 5 \\
5712 & = & 2(2100) & + & 3(390) & + & 5(60) & + & c_3(6) & & & & & & & \Rightarrow & c_3 & = & 7 \\
2100 & = & 2(602) & + & 3(180) & + & 5(50) & + & 7(12) & + & c_2(2) & & & & &  \Rightarrow & c_2 & = & 11 \\
389 & = & 2(63) & + & 3(31) & + & 5(15) & + & 7(7) & + & 11(3) & + & c_1(1) & & & \Rightarrow & c_1 & = & 13 \\
58 & = & 2(1) & + & 3(1) & + & 5(1) & + & 7(1) & + & 11(1) & + & 13(1) & + & c_0(1) & \Rightarrow & c_0 & = & 17 \\
\end{array} \]
\end{examp}

\begin{examp} \label{ex:awntni}  
In this example, we show how to extend the method to \pspl generated with integers not starting at either 0 or 1, or with input data with non-unity differentials. Since it is probably easiest to start with an integer index of 0 rather than 1 for the extension of the method, we will assign a function relating the input data to the integers $0,1,2$, etc., and we will use the values in the \awnt as in \Cref{ex:awnt}. We consider the sequence generated by the polynomial: \[3x^5 + 1x^4 + 4x^3 + 1x^2 + 5x + 9; \ x \in 3.3,3.4..3.9\] 

The difference table for this sequence and input data is given in \Cref{tbl:awntexni}, with the \md in bold. Again, the sequence values are in row 0 in accordance with \Cref{def:ddt}. This row appears directly below the integers, starting at 0, that we have calculated and assigned to the input data entries using the equation:

\begin{equation} \label{eq:gx}
g(x)=10.0(x-3.3)
\end{equation}

\noindent This allows us to use the method using the \awnt.
\vspace{0.1em}

\begin{table}[!htbp]
\centering
\caption{The Difference Table for $3x^5 + 1x^4 + 4x^3 + 1x^2 + 5x + 9; \ x \in 3.3,3.4..3.9$}  
\small                                                                
\label{tbl:awntexni}
\begin{tabular}{cccccccccccccl} 
3.3 & & \mxem 3.4 & & \mxem 3.5 & & \mxem 3.6 & & \mxem 3.7 & & \mxem 3.8 & & \mxem 3.9 & $\leftarrow x$ \\
0 & & \mxem 1 & & \mxem 2 & & \mxem 3 & & \mxem 4 & & \mxem 5 & & \mxem 6 & $\leftarrow g(x)$ \\
\textbf{1472.79189} & & \mxem 1691.47232 & & \mxem 1935.96875 & & \mxem 2208.53088 & & \mxem 2511.53681 & & \mxem 2847.49664 & & \mxem 3219.05607 \\
& \mxem \textbf{218.68043} & & \mxem 244.49643 & & \mxem 272.56213 & & \mxem 303.00593 & & \mxem 335.95983 & & \mxem 371.55943 \\
& & \mxem \textbf{25.816} & & \mxem 28.0657 & & \mxem 30.4438 & & \mxem 32.9539 & & \mxem 35.5996 \\
& & & \mxem \textbf{2.2497} & & \mxem 2.3781 & & \mxem 2.5101 & & \mxem 2.6457 \\
& & & & \mxem \textbf{0.1284} & & \mxem 0.132 & & \mxem 0.1356 \\
& & & & & \mxem \textbf{0.0036} & & \mxem 0.0036 \\
\end{tabular}
\end{table}

\vspace{0.1em} 
Row 5 of the \dt is constant, matching the degree, $d$, of the polynomial as expected. We start at $n,k=d=5$ in \Cref{tbl:awntt}, and use column $k$ to calculate $c_{n=k}$ of the polynomial using the entries in the table as multipliers for each $c_n, n \in 1..5$. By inspecting the first element in row 0 of \Cref{tbl:awntexni}, we see that $c_0$ is 1472.79189.
Furthermore, we have:
\vspace{0.2em}
\footnotesize
\[ \arraycolsep=3.0pt
\begin{array}{ccccccccccccccc} 
0.0036 & = & c_5(120) & & & & & & & & &  \Rightarrow & c_5 & = & 0.00003 \\
0.1284 & = & 0.00003(240) & + & c_4(24) & & & & & & &  \Rightarrow & c_4 & = & 0.00505 \\
2.2497 & = & 0.00003(150) & + & 0.00505(36) & + & c_3(6) & & & & &  \Rightarrow & c_3 & = & 0.3439 \\
25.816 & = & 0.00003(30) & + & 0.00505(14) & + & 0.3439(6) & + & c_2(2) & & & \Rightarrow & c_2 & = & 11.8405 \\
218.68043 & = & 0.00003(1) & + & 0.00505(1) & + & 0.3439(1) & + & 11.8405(1) & + & c_1(1) & \Rightarrow & c_1 & = & 206.49095 \\
\end{array} \]
\normalsize

\vspace{0.1em}
The value of the sequence for other input values of $x$ may be calculated using these coefficients with an input value of $g(x)$ as given in \Cref{eq:gx}. Alternatively, we could calculate the coefficients for use with $x$ directly, as opposed to $g(x)$, by symbolic evaluation of: 
\[0.00003(g(x))^5+0.00505(g(x))^4+0.3439(g(x))^3+11.8405(g(x))^2+206.49095(g(x))+1472.79189\]
which simplifies to:
\[3x^5 + 1x^4 + 4x^3 + 1x^2 + 5x + 9\]

\vspace{0.1em} 
Obviously, the same method may be used for integer generated sequences with a starting integer value other than 0 or 1 by making a substitution of $g(x)=x-y$, with $y$ as the appropriate integer. The value of $y$ will depend upon whether the \mwnt or the \awnt was used in the calculation, and upon the starting integer value of the sequence.
\end{examp}

\section{Mathematical Basis} \label{sec:mb}

\subsection{A Formula for the \mdc of the \dtc of a \psc} \label{ss:fmd}
A \ps, given by $a_i, a_{i+1},a_{i+2}$, etc., has a \dt as defined in \Cref{def:ddt}. A general example showing $a_0,a_1,a_2, \text{ and } a_3$ was given in \Cref{tbl:gdt}. It is known \cite{bbqgmds} \cite{bbcmd} that the $k^{th}$ term of the \md of the sequence of the \dt is:

\[
D_k=\sum_{i=0}^k(-1)^{k+i} \binom{k}{i} a_i \qquad k \geq 0
\]

\noindent
This may be proved via induction, or by the method given by Graham, et al., (1994) \cite{bbcmp}. If we multiply by $(-1)^{-2i}=1$, for each $i$ in turn, we get:

\begin{equation} \label{eq:md}
D_k=\sum_{i=0}^k(-1)^{k-i} \binom{k}{i} a_i \qquad k \geq 0
\end{equation}

\subsection{Using the \mdc with the \awnt} \label{ss:mdawnt}

\subsubsection{Mathematical Derivation:} \label{sss:awntmd}
In this section we derive the expressions linking the \md of the \dt of a \ps and the polynomial's coefficient multipliers to the equation for the \awnt, $\awntm$. This leads to the method of solving for the polynomial's unknown coefficients as shown in \Cref{ex:awnt}. 

First, since the starting integer used to generate the polynomial is $x=0$, it is obvious that $c_0=a_0$. Furthermore, recalling \Cref{eq:md}, the equation for the $k^{th}$ term of the \md, and substituting $a_i=\sum_{n=0}^d c_n i^n$, yields:

\[
D0_k=\sum_{i=0}^k (-1)^{k-i} \binom{k}{i} \sum_{n=0}^d c_n i^n \qquad k \geq 0
\] 

\noindent where $d$ is the degree of the polynomial.
\vspace{0.1em}

Since both sums have a finite number of terms, and due to the commutative property of multiplication and addition, and the distributive property of multiplication over addition, we may rearrange the above equation into ($c_0$ is left out as it is already known from the first value in the main diagonal):

\begin{equation} \label{eq:awntmde}
D0_k=\sum_{n=1}^d c_n \ \underbrace{\sum_{i=0}^k(-1)^{k-i} \binom{k}{i} i^n}_{\awntm} \qquad n \geq 1, \ k \geq 1
\end{equation}

This equation shows that that the $k^{th}$ main diagonal element value for $k  \geq 1$, is composed of $c_n, \ n \in 1..d$, multiplied by $\sum_{i=0}^k (-1)^{k-i} \binom{k}{i} i^n$. This is column $k$, with corresponding row, $n$, in \awnta, and is seen by comparison of \Cref{eq:awntmde} with \Cref{eq:awnt}. 

\subsubsection{Solution Procedure using \boldmath \texorpdfstring{\awntm}{AWNT(n,k)} \unboldmath \hspace{-0.3em}, \awnta:} \label{sss:awntsp}
Therefore, to find the unknown coefficients of a polynomial sequence generated with integers starting at 0, first construct the \dt until the elements of a row are all constant. The degree, $d$, of the polynomial is the row number, as defined in \Cref{def:ddt}, of the difference triangle with the constant values. Then start with column $k=d$ and row $n=d$ in \awnta, and for each $c_n, n  \in 1..d$, assign the value of \awntm as a multiplier to $c_n$ (moving up the column is perhaps easiest) and equate it to the main diagonal value for row $k=n=d$ in the \dt. \awntm[n][d] will have multipliers of 0 for all $c_n$ except for $c_d$, allowing for easy calculation of $c_d$.

Next, move back to column $k=d-1$ and starting at row $n=d$ assign the value of multipliers to each $c_n$, again moving up the column. Since the value of $c_d$ is known, equating the coefficients and multipliers to the main diagonal value in row $d-1$ will leave $c_{d-1}$ as the only unknown. Continue this process backwards to column $k=1$ to solve for the coefficients down to $c_1$. The value of $c_0$ is equal to the value in row 0 of the \dt (the first term of the sequence), and the solution is complete. See \Cref{ex:awnt} for a worked example using this procedure.

\subsection{Using the \mdc with the \mwnt} \label{ss:mdmwnt}

\subsubsection{Mathematical Derivation:} \label{sss:mwntmd}
In this section we derive the expressions linking the \md of the \dt of a \ps and the polynomial's coefficient multipliers to the equation for the \mwnt, $\mwntm$. This leads to the method of solving for the polynomial's unknown coefficients as shown in \Cref{ex:mwnt}. In this case, the starting integer used to generate the polynomial is $x=1$, and we conveniently refer to the terms of the sequence as $a_1, a_2, a_3$, etc. If we look at the $m^{th}$ term of the \md from \Cref{eq:md}, we get:

\[
D1_m=\sum_{j=0}^m(-1)^{m-j} \binom{m}{j}a_{j+1} \qquad m \geq 0
\]

\noindent where:

\[
a_{j+1}=\sum_{q=0}^d c_q (j+1)^q
\]

\noindent and $d$ is the degree of the polynomial, as before. Substituting, we get:

\[
D1_m=\sum_{j=0}^m(-1)^{m-j} \binom{m}{j} \sum_{q=0}^d c_q (j+1)^q \qquad m \geq 0
\]

Since both sums have a finite number of terms, and due to the commutative property of multiplication and addition, and the distributive property of multiplication over addition, we may rearrange the above equation into:

\begin{align} \label{eq:intermed}
D1_m &= \sum_{q=0}^d c_q \sum_{j=0}^m(-1)^{m-j} \binom{m}{j} (j+1)^q & \text{Let } i &=j+1 \Rightarrow j=i-1: \nonumber \\
D1_m &= \sum_{q=0}^d c_q \sum_{i=1}^{m+1}(-1)^{m-i+1} \binom{m}{i-1} i^q & \text{Let } k&=m+1 \Rightarrow m=k-1: \nonumber \\ 
D1_{k-1} &= \sum_{q=0}^d c_q \sum_{i=1}^{k}(-1)^{k-i} \binom{k-1}{i-1} i^q & k \geq 1
\end{align}

We now need to show that the following relationship involving the right hand sum of \Cref{eq:intermed} is valid:

\begin{equation} \label{eq:qpone}
\sum_{i=1}^{k}(-1)^{k-i} \binom{k-1}{i-1} i^q = \frac{1}{k} \sum_{i=0}^{k} (-1)^{k-i} \binom{k}{i} i^{q+1} \qquad q \geq 0, \ k \geq 1
\end{equation}

First, on the right hand side, the $i=0$ term is 0 since $0^{q+1}=0$. The rest of the terms, with $i \geq 1$, are equal on a term by term basis, shown as follows:

\begin{align*}
(-1)^{k-i} \binom{k-1}{i-1}i^q & \overset{?}{=} \frac{1}{k} (-1)^{k-i} \binom{k}{i} i^{q+1}  \\  
\binom{k-1}{i-1}i^q & \overset{?}{=} \frac{1}{k} \binom{k}{i} i^{q+1} \\ 
\frac{(k-1)! \ i^q}{(k-1-(i-1))! \ (i-1)!} & \overset{?}{=} \frac{1}{k} \ \frac{k! \ i^{q+1}}{(k-i)! \ i!} \\ 
\frac{(k-1)! \ i^q}{(k-i)! \ (i-1)!} & \overset{?}{=} \frac{(k-1)! \ i^{q+1}}{(k-i)! \ i!} \\
\frac{i^q}{(i-1)!} & \overset{?}{=} \frac{i^q \ i}{i!} \\
\frac{i^q}{(i-1)!} & \overset{\checkmark}{=} \frac{i^q}{(i-1)!} 
\end{align*}

\noindent which confirms the relationship (again, with $i \geq 1$). We now substitute the right side expression of \Cref{eq:qpone} for the right side sum of \Cref{eq:intermed} and get:

\begin{align} \label{eq:mwntmde}
D1_{k-1} &= \sum_{q=0}^d c_q \ \frac{1}{k} \sum_{i=0}^{k} (-1)^{k-i} \binom{k}{i} i^{q+1} & q \geq 0, \ k \geq 1 & & \text{Let } n=q+1 \Rightarrow q=n-1 \nonumber \\
D1_{k-1} &= \sum_{n=1}^{d+1} c_{n-1} \ \underbrace{\frac{1}{k} \sum_{i=0}^{k} (-1)^{k-i} \binom{k}{i} i^n}_{\mwntm} & n \geq 1, \ k \geq 1 &
\end{align}

This equation shows that for $k \geq 1$, the $(k-1)^{th}$ main diagonal element value is composed of $c_{n-1}, \ n \in 1..(d+1)$, multiplied by $\frac{1}{k} \sum_{i=0}^{k} (-1)^{(k-i)} \binom{k}{i} i^n$. This is column, $k$, with corresponding row, $n$, in \mwnta, as seen by comparison of \Cref{eq:mwntmde} with \Cref{eq:mwnt}.

\subsubsection{Solution Procedure using \boldmath \texorpdfstring{\mwntm}{MWNT(n,k)} \unboldmath \hspace{-0.3em}, \mwnta:} \label{sss:mwntsp}
Therefore, to find the unknown coefficients of a polynomial sequence generated with integers starting at 1, first construct the \dt until the elements of a row are all constant. The degree, $d$, of the polynomial is the row number, as defined in \Cref{def:ddt}, of the difference triangle with the constant values. Then start with column $k=d+1$ and row $n=d+1$ in \mwnta, and for each $c_{n-1}, n  \in 1..(d+1)$, assign the value of \mwntm as a multiplier to $c_{n-1}$ (moving up the column) and equate it to the main diagonal value for row $k-1=n-1=d$ in the \dt. \mwntm[n][d+1] will have multipliers of 0 for all $c_{n-1}$ except for $c_d$, allowing for easy calculation of $c_d$. 

Next, move back to column $k=d$ and starting at row $n=d+1$ assign the value of multipliers to each $c_{n-1}$, again moving up the column. Since the value of $c_d$ is known, equating the coefficients and multipliers to the main diagonal value in row $d-1$ will leave $c_{d-1}$ as the only unknown. Continue this process backwards to column $k=1$ to solve for the coefficients down to $c_0$, and the solution is complete. See \Cref{ex:mwnt} for a worked example using this procedure.

\subsection{Gaining Insight Using \efdt} \label{ss:mbefdt}
\efdt as presented by \gaq \cite{bbqgefdt} states that given $f(x)=\sum_{j=0}^d c_j x^j$ then:

\[
\sum_{i=0}^k(-1)^i \binom{k}{i}f(i) = 
\begin{cases}
0, & 0 \leq d < k \\
(-1)^k k! \ c_k, & d=k
\end{cases}
\]

\vspace{0.2em}
\noindent The authors use \efdt and let: \[f(x)=(z-bx)^n \ \Rightarrow \ f(i)=(z-bi)^n, \ n \in \mathbb{Z}_{\geq 0}\] to derive the following equation \cite{bbqgmme}: 

\begin{equation} \label{eq:sec}
\sum_{i=0}^k (-1)^i \binom{k}{i} (z-bi)^n = 
\begin{cases}
0, & n < k \\
b^k \ k!, & n=k
\end{cases}
\end{equation}

\vspace{0.2em}
\noindent By setting $b=-1$ and $z=0$ we can derive \Cref{eq:main} below as follows: 

\[
\sum_{i=0}^k (-1)^i \binom{k}{i} i^n =
\begin{cases}
0, & n < k \\
(-1)^k k!, & n=k
\end{cases}
\]

\noindent If we multiply each side by $(-1)^k$, we get:

\[
\sum_{i=0}^k (-1)^{k+i} \binom{k}{i} i^n =
\begin{cases}
0, & n < k \\
(-1)^{2k} k!, & n=k
\end{cases}
\]

\noindent If we multiply the left side by $(-1)^{-2i}=1$, for each $i$ in turn, and since $(-1)^{2k}=1$, we get:

\begin{equation} \label{eq:main}
\sum_{i=0}^k (-1)^{k-i} \binom{k}{i} i^n =
\begin{cases}
0, & n < k \\
k!, & n=k
\end{cases}
\end{equation}

We can use \Cref{eq:sec} and \Cref{eq:main} to gain insight into the structure of the contents of the \mwnt (\Cref{tbl:mwntt}), the \awnt (\Cref{tbl:awntt}) and the \dt of a sequence. \Cref{eq:main} shows why the factorials of numbers appear on the right diagonal of both triangles, starting from 0! in \Cref{tbl:mwntt}, and from 1! in \Cref{tbl:awntt}. It also explains the zeros in the tables when $n<k$, and the rows of 0 that one would get by continuing the \dt past the constant row, as $n<k$.

\Cref{eq:sec} may be used to explain the constant row of the \dt. In that equation, $z$ may be taken to be any number (not just 0), so along with having $b=-1$, this explains why successive terms in the row are all equal (and related to the factorial), and that the method of isolating the coefficients could be used with any row of the \dt because the coefficient multipliers would also remain 0 for $n<k$.\footnote{A multiplication by $(-1)^k$ is assumed as was done in proceeding from \Cref{eq:sec} to \Cref{eq:main} in order to match the form of the diagonal terms where the multipliers of the $i=k$ term are positive (so that the factorials are all positive).} However, different tables of numbers would need to be used for the multipliers (the factorials and zeros would still be in place), or the multipliers could just be calculated from the terms in \Cref{eq:sec} with the appropriate values for $z$ and $b$ -- see \Cref{sec:eor} for the terms used in \Cref{ex:awnt} with \awnta ($z=0, \ b=-1$).

\section{\texorpdfstring{\Cref{ex:awnt}}{Example 1} Revisited} \label{sec:eor}
In this section, we will show \Cref{ex:awnt} in full, per \Cref{eq:awntmde} with all terms shown. Along with \Cref{eq:main}, this will hopefully provide a more complete view of how the method works. So we have, with the binomial coefficients shown in bold:

\arraycolsep=0.5pt 
\[ \begin{array}{ccccccccccccccccccccl}    
2880 & = & c_6 & \cdot & (\mathbf{1} \cdot 6^6 & - & \mathbf{6} \cdot 5^6 & + & \mathbf{15} \cdot 4^6 & - & \mathbf{20} \cdot 3^6 & + & \mathbf{15} \cdot 2^6 & - & \mathbf{6} \cdot 1^6 & + & \mathbf{1} \cdot 0^6) & = & c_6 \cdot \awntm[6][6] & =  & c_6 \cdot 720 \\
& + & c_5 & \cdot & (\mathbf{1} \cdot 6^5 & - & \mathbf{6} \cdot 5^5 & + & \mathbf{15} \cdot 4^5 & - & \mathbf{20} \cdot 3^5 & + & \mathbf{15} \cdot 2^5 & - & \mathbf{6} \cdot 1^5 & + & \mathbf{1} \cdot 0^5) & = & c_5 \cdot \awntm[5][6] & =  & c_5 \cdot 0 \\
& + & c_4 & \cdot & (\mathbf{1} \cdot 6^4 & - & \mathbf{6} \cdot 5^4 & + & \mathbf{15} \cdot 4^4 & - & \mathbf{20} \cdot 3^4 & + & \mathbf{15} \cdot 2^4 & - & \mathbf{6} \cdot 1^4 & + & \mathbf{1} \cdot 0^4) & = & c_4 \cdot \awntm[4][6] & =  & c_4 \cdot 0 \\
& + & c_3 & \cdot & (\mathbf{1} \cdot 6^3 & - & \mathbf{6} \cdot 5^3 & + & \mathbf{15} \cdot 4^3 & - & \mathbf{20} \cdot 3^3 & + & \mathbf{15} \cdot 2^3 & - & \mathbf{6} \cdot 1^3 & + & \mathbf{1} \cdot 0^3) & = & c_3 \cdot \awntm[3][6] & =  & c_3 \cdot 0 \\
& + & c_2 & \cdot & (\mathbf{1} \cdot 6^2 & - & \mathbf{6} \cdot 5^2 & + & \mathbf{15} \cdot 4^2 & - & \mathbf{20} \cdot 3^2 & + & \mathbf{15} \cdot 2^2 & - & \mathbf{6} \cdot 1^2 & + & \mathbf{1} \cdot 0^2) & = & c_2 \cdot \awntm[2][6] & =  & c_2 \cdot 0 \\
& + & c_1 & \cdot & (\mathbf{1} \cdot 6^1 & - & \mathbf{6} \cdot 5^1 & + & \mathbf{15} \cdot 4^1 & - & \mathbf{20} \cdot 3^1 & + & \mathbf{15} \cdot 2^1 & - & \mathbf{6} \cdot 1^1 & + & \mathbf{1} \cdot 0^1) & = & c_1 \cdot \awntm[1][6] & =  & c_1 \cdot 0 \\
c_6 \Rightarrow 4 \\
\\
	7800 & = & 4 & \cdot & (\mathbf{1} \cdot 5^6 & - & \mathbf{5} \cdot 4^6 & + & \mathbf{10} \cdot 3^6 & - & \mathbf{10} \cdot 2^6 & + & \mathbf{5} \cdot 1^6 & - & \mathbf{1} \cdot 0^6) & & & = & 4 \cdot \awntm[6][5] & =  & 4 \cdot 1800 \\
& + & c_5 & \cdot & (\mathbf{1} \cdot 5^5 & - & \mathbf{5} \cdot 4^5 & + & \mathbf{10} \cdot 3^5 & - & \mathbf{10} \cdot 2^5 & + & \mathbf{5} \cdot 1^5 & - & \mathbf{1} \cdot 0^5) & & & = & c_5 \cdot \awntm[5][5] & =  & c_5 \cdot 120 \\
& + & c_4 & \cdot & (\mathbf{1} \cdot 5^4 & - & \mathbf{5} \cdot 4^4 & + & \mathbf{10} \cdot 3^4 & - & \mathbf{10} \cdot 2^4 & + & \mathbf{5} \cdot 1^4 & - & \mathbf{1} \cdot 0^4) & & & = & c_4 \cdot \awntm[4][5] & =  & c_4 \cdot 0 \\
& + & c_3 & \cdot & (\mathbf{1} \cdot 5^3 & - & \mathbf{5} \cdot 4^3 & + & \mathbf{10} \cdot 3^3 & - & \mathbf{10} \cdot 2^3 & + & \mathbf{5} \cdot 1^3 & - & \mathbf{1} \cdot 0^3) & & & = & c_3 \cdot \awntm[3][5] & =  & c_3 \cdot 0 \\
& + & c_2 & \cdot & (\mathbf{1} \cdot 5^2 & - & \mathbf{5} \cdot 4^2 & + & \mathbf{10} \cdot 3^2 & - & \mathbf{10} \cdot 2^2 & + & \mathbf{5} \cdot 1^2 & - & \mathbf{1} \cdot 0^2) & & & = & c_2 \cdot \awntm[2][5] & =  & c_2 \cdot 0 \\
& + & c_1 & \cdot & (\mathbf{1} \cdot 5^1 & - & \mathbf{5} \cdot 4^1 & + & \mathbf{10} \cdot 3^1 & - & \mathbf{10} \cdot 2^1 & + & \mathbf{5} \cdot 1^1 & - & \mathbf{1} \cdot 0^1) & & & = & c_1 \cdot \awntm[1][5] & =  & c_1 \cdot 0 \\
c_5 \Rightarrow 5 \\
\\

7584 & = & 4 & \cdot & (\mathbf{1} \cdot 4^6 & - & \mathbf{4} \cdot 3^6 & + & \mathbf{6} \cdot 2^6 & - & \mathbf{4} \cdot 1^6 & + & \mathbf{1} \cdot 0^6) & & & & & = & 4 \cdot \awntm[6][4] & =  & 4 \cdot 1560 \\
& + & 5 & \cdot & (\mathbf{1} \cdot 4^5 & - & \mathbf{4} \cdot 3^5 & + & \mathbf{6} \cdot 2^5 & - & \mathbf{4} \cdot 1^5 & + & \mathbf{1} \cdot 0^5) & & & & & = & 5 \cdot \awntm[5][4] & =  & 5 \cdot 240 \\
& + & c_4 & \cdot & (\mathbf{1} \cdot 4^4 & - & \mathbf{4} \cdot 3^4 & + & \mathbf{6} \cdot 2^4 & - & \mathbf{4} \cdot 1^4 & + & \mathbf{1} \cdot 0^4) & & & & & = & c_4 \cdot \awntm[4][4] & =  & c_4 \cdot 24 \\
& + & c_3 & \cdot & (\mathbf{1} \cdot 4^3 & - & \mathbf{4} \cdot 3^3 & + & \mathbf{6} \cdot 2^3 & - & \mathbf{4} \cdot 1^3 & + & \mathbf{1} \cdot 0^3) & & & & & = & c_3 \cdot \awntm[3][4] & =  & c_3 \cdot 0 \\
& + & c_2 & \cdot & (\mathbf{1} \cdot 4^2 & - & \mathbf{4} \cdot 3^2 & + & \mathbf{6} \cdot 2^2 & - & \mathbf{4} \cdot 1^2 & + & \mathbf{1} \cdot 0^2) & & & & & = & c_2 \cdot \awntm[2][4] & =  & c_2 \cdot 0 \\
& + & c_1 & \cdot & (\mathbf{1} \cdot 4^1 & - & \mathbf{4} \cdot 3^1 & + & \mathbf{6} \cdot 2^1 & - & \mathbf{4} \cdot 1^1 & + & \mathbf{1} \cdot 0^1) & & & & & = & c_1 \cdot \awntm[1][4] & =  & c_1 \cdot 0 \\
c_4 \Rightarrow 6 \\
\\

3168 & = & 4 & \cdot & (\mathbf{1} \cdot 3^6 & - & \mathbf{3} \cdot 2^6 & + & \mathbf{3} \cdot 1^6 & - & \mathbf{1} \cdot 0^6) & & & & & & & = & 4 \cdot \awntm[6][3] & =  & 4 \cdot 540 \\
& + & 5 & \cdot & (\mathbf{1} \cdot 3^5 & - & \mathbf{3} \cdot 2^5 & + & \mathbf{3} \cdot 1^5 & - & \mathbf{1} \cdot 0^5) & & & & & & & = & 5 \cdot \awntm[5][3] & =  & 5 \cdot 150 \\
& + & 6 & \cdot & (\mathbf{1} \cdot 3^4 & - & \mathbf{3} \cdot 2^4 & + & \mathbf{3} \cdot 1^4 & - & \mathbf{1} \cdot 0^4) & & & & & & & = & 6 \cdot \awntm[4][3] & =  & 6 \cdot 36 \\
& + & c_3 & \cdot & (\mathbf{1} \cdot 3^3 & - & \mathbf{3} \cdot 2^3 & + & \mathbf{3} \cdot 1^3 & - & \mathbf{1} \cdot 0^3) & & & & & & & = & c_3 \cdot \awntm[3][3] & =  & c_3 \cdot 6 \\
& + & c_2 & \cdot & (\mathbf{1} \cdot 3^2 & - & \mathbf{3} \cdot 2^2 & + & \mathbf{3} \cdot 1^2 & - & \mathbf{1} \cdot 0^2) & & & & & & & = & c_2 \cdot \awntm[2][3] & =  & c_2 \cdot 0 \\
& + & c_1 & \cdot & (\mathbf{1} \cdot 3^1 & - & \mathbf{3} \cdot 2^1 & + & \mathbf{3} \cdot 1^1 & - & \mathbf{1} \cdot 0^1) & & & & & & & = & c_1 \cdot \awntm[1][3] & =  & c_1 \cdot 0 \\
c_3 \Rightarrow 7 \\
\\

540 & = & 4 & \cdot & (\mathbf{1} \cdot 2^6 & - & \mathbf{2} \cdot 1^6 & + & \mathbf{1} \cdot 0^6) & & & & & & & & & = & 4 \cdot \awntm[6][2] & =  & 4 \cdot 62 \\
& + & 5 & \cdot & (\mathbf{1} \cdot 2^5 & - & \mathbf{2} \cdot 1^5 & + & \mathbf{1} \cdot 0^5) & & & & & & & & & = & 5 \cdot \awntm[5][2] & =  & 5 \cdot 30 \\
& + & 6 & \cdot & (\mathbf{1} \cdot 2^4 & - & \mathbf{2} \cdot 1^4 & + & \mathbf{1} \cdot 0^4) & & & & & & & & & = & 6 \cdot \awntm[4][2] & =  & 6 \cdot 14 \\
& + & 7 & \cdot & (\mathbf{1} \cdot 2^3 & - & \mathbf{2} \cdot 1^3 & + & \mathbf{1} \cdot 0^3) & & & & & & & & & = & 7 \cdot \awntm[3][2] & =  & 7 \cdot 6 \\
& + & c_2 & \cdot & (\mathbf{1} \cdot 2^2 & - & \mathbf{2} \cdot 1^2 & + & \mathbf{1} \cdot 0^2) & & & & & & & & & = & c_2 \cdot \awntm[2][2] & =  & c_2 \cdot 2 \\
& + & c_1 & \cdot & (\mathbf{1} \cdot 2^1 & - & \mathbf{2} \cdot 1^1 & + & \mathbf{1} \cdot 0^1) & & & & & & & & & = & c_1 \cdot \awntm[1][2] & =  & c_1 \cdot 0 \\
c_2 \Rightarrow 8 \\
\\

39 & = & 4 & \cdot & (\mathbf{1} \cdot 1^6 & - & \mathbf{1} \cdot 0^6) & & & & & & & & & & & = & 4 \cdot \awntm[6][1] & =  & 4 \cdot 1 \\
& + & 5 & \cdot & (\mathbf{1} \cdot 1^5 & - & \mathbf{1} \cdot 0^5) & & & & & & & & & & & = & 5 \cdot \awntm[5][1] & =  & 5 \cdot 1 \\
& + & 6 & \cdot & (\mathbf{1} \cdot 1^4 & - & \mathbf{1} \cdot 0^4) & & & & & & & & & & & = & 6 \cdot \awntm[4][1] & =  & 6 \cdot 1 \\
& + & 7 & \cdot & (\mathbf{1} \cdot 1^3 & - & \mathbf{1} \cdot 0^3) & & & & & & & & & & & = & 7 \cdot \awntm[3][1] & =  & 7 \cdot 1 \\
& + & 8 & \cdot & (\mathbf{1} \cdot 1^2 & - & \mathbf{1} \cdot 0^2) & & & & & & & & & & & = & 8 \cdot \awntm[2][1] & =  & 8 \cdot 1 \\
& + & c_1 & \cdot & (\mathbf{1} \cdot 1^1 & - & \mathbf{1} \cdot 0^1) & & & & & & & & & & & = & c_1 \cdot \awntm[1][1] & =  & c_1 \cdot 1 \\
c_1 \Rightarrow 9
\end{array} \]

\noindent We know that $c_0=10$ for the reasons stated before in \Cref{ex:awnt}. This completes the solution. 

\vspace{1em}
Further notes: Obviously, $c_0$ is absent from any row number greater than 0 in the \dt as it is subtracted out from each term in row 1. It should also be noted that $c_1$ will be absent from any row number greater than 1 in the \dt because in row 1 there will be a difference of $c_1 \cdot 1x$, which will subtracted out from each term in row 2. This is not the case for the higher order terms.

\section{\texorpdfstring{\Cref{ex:mwnt}}{Example 2} (Partially) Revisited} \label{sec:etr}
In this section, we will show a partial working of \Cref{ex:mwnt}. Although only partial, this presentation of the solution will hopefully clarify how the method works by going into more detail in the areas that are covered, especially when compared to the full solution to \Cref{ex:awnt} given above. Our focus will be on \Cref{eq:intermed} and \Cref{eq:mwntmde}, for the calculation of $c_6$.

\Cref{eq:intermed} gives the expression derived for the \md values and the coefficients per \Cref{eq:md}, given that the starting integer used to generate the sequence is $x=1$. Since $k=7$, and with $q$ taken appropriately for each coefficient, we have, (partially):

\arraycolsep=1pt 
\[ \begin{array}{ccccccccccccccccccl}    
1440 & = & c_6 & \cdot & (\mathbf{1} \cdot 7^6 & - & \mathbf{6} \cdot 6^6 & + & \mathbf{15} \cdot 5^6 & - & \mathbf{20} \cdot 4^6 & + & \mathbf{15} \cdot 3^6 & - & \mathbf{6} \cdot 2^6 & + & \mathbf{1} \cdot 1^6) & =  & c_6 \cdot 720 \\
& + & c_5 & \cdot & (\mathbf{1} \cdot 7^5 & - & \mathbf{6} \cdot 6^5 & + & \mathbf{15} \cdot 5^5 & - & \mathbf{20} \cdot 4^5 & + & \mathbf{15} \cdot 3^5 & - & \mathbf{6} \cdot 2^5 & + & \mathbf{1} \cdot 1^5) & =  & c_5 \cdot 0 \\
& \vdots
\end{array} \]

\noindent \Cref{eq:mwntmde}, with $n=q+1$ taken appropriately for each coefficient, gives (partially):
\arraycolsep=0.0pt 
\[ \begin{array}{ccccccccccccccccccccccl}    
1440 & = & c_6 & \cdot \frac{1}{7} & (\mathbf{1} \cdot 7^7 & - & \mathbf{7} \cdot 6^7 & + & \mathbf{21} \cdot 5^7 & - & \mathbf{35} \cdot 4^7 & + & \mathbf{35} \cdot 3^7 & - & \mathbf{21} \cdot 2^7 & + & \mathbf{7} \cdot 1^7 & - & \mathbf{1} \cdot 0^7) & = & c_6 \cdot \mwntm[7][7] & =  & c_6 \cdot 720 \vspace{0.2em} \\
& + & c_5 &\cdot \frac{1}{7} & (\mathbf{1} \cdot 7^6 & - & \mathbf{7} \cdot 6^6 & + & \mathbf{21} \cdot 5^6 & - & \mathbf{35} \cdot 4^6 & + & \mathbf{35} \cdot 3^6 & - & \mathbf{21} \cdot 2^6 & + & \mathbf{7} \cdot 1^6 & - & \mathbf{1} \cdot 0^6) & = & c_5 \cdot \mwntm[6][7] & =  & c_5 \cdot 0 \\
& \vdots
\end{array} \]

\section{Closing Remarks} \label{sec:cr}
We have presented a method of solving for the unknown coefficients of a \ps using two of the \wnts and the \md of the sequence's \dt. However, we are unsure if our description of the method in this paper meets the rigorous requirements for a mathematical proof.\footnote{Our experience is in engineering, not in producing ironclad mathematical proofs. The recognition of the method as presented in this paper came about upon investigation of the results of the statistical analysis of an electronics manufacturing process.} Regarding the use of the method described in \Cref{ss:mdawnt} and \Cref{ss:mdmwnt}, if we have not met the requirements, we feel that we are fairly close in doing so. We are quite sure that we have not met the requirements of mathematical rigor regarding extending the method as in \Cref{ex:awntni}. However, we believe that the extension of the method is valid and that any result obtained in practice may be checked for validity on a case by case basis.

We welcome papers that address any shortcomings in this paper, with due credit going to the authors.\footnote{We are assuming, perhaps incorrectly, that members of the mathematical community with the necessary skills and experience to write such papers will also deem it worthwhile to do so.}\footnote{We also note that while we are interested in real numbers, we surmise that mathematicians may wish to extend the method to include complex numbers, if possible.} Certainly, we feel that the \wnts are worthy of wider recognition and perhaps of standard definition and notation as well.

\clearpage 

\def\bibindent{1em}

\end{document}